\documentclass[12pt,oneside,english]{amsart}
\usepackage[T1]{fontenc}
\usepackage[latin1]{inputenc}
\usepackage{geometry}
\usepackage{setspace}
\usepackage{amssymb}

\makeatletter
 \theoremstyle{plain}
\newtheorem{thm}{Theorem}[section]
  \theoremstyle{plain}
  \newtheorem{lem}[thm]{Lemma}
  \theoremstyle{plain}
  \newtheorem{prop}[thm]{Proposition}
  \theoremstyle{plain}
  \newtheorem{cor}[thm]{Corollary}
 \theoremstyle{definition}
 \newtheorem*{defn*}{Definition}
  \theoremstyle{remark}
  \newtheorem*{rem*}{Remark}

\numberwithin{equation}{section}

\usepackage{babel}
\makeatother
\begin{document}

\title{limit laws for boolean convolutions}

\author{jiun-chau wang}

\date{July 28th, 2007}

\begin{abstract}
We study the distributional behavior for products and sums of boolean
independent random variables in a general infinitesimal triangular
array. We show that the limit laws of boolean convolutions are determined
by the limit laws of free convolutions, and vice versa. We further
use these results to demonstrate several connections between the limiting
distributional behavior of classical convolutions and that of boolean
convolutions. The proof of our results is based on the analytical
apparatus developed in \cite{BerJcMul,BerJcAdd} for free convolutions. 
\end{abstract}
\maketitle

\section{Introduction}

\def\Utimes{\cup\kern-.50em\lower-.4ex\hbox{$_\times$}}\def\utimes{\cup\kern-.86em\lower-.7ex\hbox{$_\times$}\,}\def\Utime{\cup\!\!\!\!\lower-.6ex\hbox{$_\times$}\,}Denote
by $\mathcal{M}_{\mathbb{R}}$ the collection of all Borel probability
measures on the real line $\mathbb{R}$, and by $\mathcal{M}_{\mathbb{T}}$
Borel probability measures on the unit circle $\mathbb{T}$. The classical
convolution $*$ for elements in $\mathcal{M}_{\mathbb{R}}$ corresponds
to the addition of independent real random variables, and the convolution
$\circledast$ for measures in $\mathcal{M}_{\mathbb{T}}$ corresponds
to the multiplication of independent circle-valued random variables.
A binary operation $\uplus$ on $\mathcal{M}_{\mathbb{R}}$, called
\emph{additive boolean convolution}, was introduced by Speicher and
Woroudi \cite{SpeWoro}. They also showed that it corresponds to the
addition of random variables belonging to algebras which are boolean
independent. Later Franz \cite{Franz} introduced the concept of \emph{multiplicative
boolean convolution} $\Utime$ for measures in $\mathcal{M}_{\mathbb{T}}$,
which is a multiplicative counterpart of the additive boolean convolution.
As shown by Voiculescu \cite{Voiadd,Voimul}, there are two other
convolutions defined respectively for measures on $\mathbb{R}$ and
$\mathbb{T}$. These are additive free convolution $\boxplus$ and
multiplicative free convolution $\boxtimes$. 

The main purpose of this paper is to investigate the limiting distributional
behavior for boolean convolutions of measures in an infinitesimal
triangular array. Let $\{ k_{n}\}_{n=1}^{\infty}$ be a sequence of
natural numbers. A triangular array $\{\mu_{nk}:n\in\mathbb{N},1\le k\leq k_{n}\}\subset\mathcal{M}_{\mathbb{T}}$
is said to be \emph{infinitesimal} if \[
\lim_{n\rightarrow\infty}\max_{1\leq k\leq k_{n}}\mu_{nk}(\{\zeta\in\mathbb{T}:\,\left|\zeta-1\right|\geq\varepsilon\})=0,\]
for every $\varepsilon>0$. Given such an array and a sequence $\{\lambda_{n}\}_{n=1}^{\infty}\subset\mathbb{T}$,
define \[
\mu_{n}=\delta_{\lambda_{n}}\utimes\mu_{n1}\utimes\mu_{n2}\utimes\cdots\utimes\mu_{nk_{n}},\qquad\nu_{n}=\delta_{\lambda_{n}}\boxtimes\mu_{n1}\boxtimes\mu_{n2}\boxtimes\cdots\boxtimes\mu_{nk_{n}},\]
and \[
\sigma_{n}=\delta_{\lambda_{n}}\circledast\mu_{n1}\circledast\mu_{n2}\circledast\cdots\circledast\mu_{nk_{n}},\qquad n\in\mathbb{N},\]
where $\delta_{\lambda_{n}}$ is the point mass at $\lambda_{n}$.
We first prove in this paper that any weak limit of such a sequence
$\{\mu_{n}\}_{n=1}^{\infty}$ is an infinitely divisible measure.
This result may be viewed as the multiplicative boolean analogue of
Hin\v{c}in's classical theorem \cite{Hincin}. Note that the same
result for $\boxplus$ (resp., $\boxtimes$) has been proved in \cite{BerPataAddHincin}
(resp., \cite{SerBerMulHincin}). Next, we find necessary and sufficient
conditions for the weak convergence of $\mu_{n}$ to a given infinitely
divisible measure. In particular, our results show that the sequence
$\mu_{n}$ converges weakly if and only if the sequence $\nu_{n}$
converges weakly. As an application, we show that the measures $\sigma_{n}$
have a weak limit if the measures $\mu_{n}$ (or $\nu_{n}$) have
a weak limit whose first moment is not zero. Moreover, the classical
limits and the boolean limits are related in an explicit manner. We
also introduce the notion of boolean normal distributions on $\mathbb{T}$,
and we show that the sequence $\mu_{n}$ converges weakly to such
a distribution if and only if the sequence $\sigma_{n}$ converges
weakly to a normal distribution (which is the push-forward measure
of a Gaussian law on $\mathbb{R}$ via the natural homomorphism from
$\mathbb{R}$ into $\mathbb{T}$.) 

The additive version of our results were studied earlier by Bercovici
and Pata in \cite{BerPataStable} for arrays with identically distributed
rows. Thus, consider an infinitesimal array $\{\nu_{nk}\}_{n,k}\subset\mathcal{M}_{\mathbb{R}}$
with $\eta_{n}=\nu_{n1}=\nu_{n2}=\cdots=\nu_{nk_{n}}$, $n\in\mathbb{N}$.
The infinitesimality here means that \[
\lim_{n\rightarrow\infty}\max_{1\leq k\leq k_{n}}\nu_{nk}(\{ t\in\mathbb{R}:\,\left|t\right|\geq\varepsilon\})=0,\]
for every $\varepsilon>0$. Set\[
\rho_{n}=\underbrace{\eta_{n}*\eta_{n}*\cdots*\eta_{n}}_{k_{n}\:\text{times}},\qquad\tau_{n}=\underbrace{\eta_{n}\boxplus\eta_{n}\boxplus\cdots\boxplus\eta_{n}}_{k_{n}\:\text{times}}\]
and\[
\omega_{n}=\underbrace{\eta_{n}\uplus\eta_{n}\uplus\cdots\uplus\eta_{n}}_{k_{n}\:\text{times}},\qquad n\in\mathbb{N}.\]
The main result in \cite[Theorem 6.3]{BerPataStable} is the equivalences
of weak convergence among the sequences $\rho_{n}$, $\tau_{n}$ and
$\omega_{n}$. The result concerning $\rho_{n}$ and $\tau_{n}$ was
first extended to an arbitrary infinitesimal array by Chistyakov and
G\"{o}tze \cite{CG2} (see also \cite{BerJcAdd} for a different
argument.) In the last part of this paper, we show how to extend the
result regarding $\tau_{n}$ and $\omega_{n}$ to an arbitrary infinitesimal
array using the methods in \cite{BerJcAdd}. 

The remainder of this paper is organized as follows. In Section 2,
we review the analytic tools needed for the calculation of boolean
convolutions. We also describe the analytic characterization of infinite
divisibility related to the various convolutions. In Section 3 we
prove the limit theorems for arrays on $\mathbb{T}$. The results
regarding the classical convolution $\circledast$ are proved in Section
4. Finally, we present the analogous results for arrays on $\mathbb{R}$
in Section 5.

\section{Preliminaries}

The analytic methods needed for the calculation of free convolutions
was discovered by Voiculescu \cite{Voiadd,Voimul}. Likewise, the
additive boolean convolution formula was found by Speicher and Woroudi
\cite{SpeWoro}, and the basic analysis of the multiplicative boolean
convolution was done by Franz \cite{Franz} (see also the paper of
Bercovici \cite{BerBoolean} for a different approach to the calculation
of both boolean convolutions.) The details are as follows.

\subsection{Multiplicative boolean and free convolutions on the unit circle}

Denote by $\mathbb{D}$ the open unit disk of the complex plane $\mathbb{C}$,
and by $\overline{\mathbb{D}}$ the closed unit disk of $\mathbb{C}$.
For a probability measure $\mu$ supported on $\mathbb{T}$, one defines
the analytic function $B_{\mu}:\,\mathbb{D}\rightarrow\mathbb{C}$
by \[
B_{\mu}(z)=\frac{1}{z}\frac{\psi_{\mu}(z)}{1+\psi_{\mu}(z)},\qquad z\in\mathbb{D},\]
where the formula of $\psi_{\mu}$ is given by\[
\psi_{\mu}(z)=\int_{\mathbb{T}}\frac{\zeta z}{1-\zeta z}\, d\mu(\zeta).\]
Note that \begin{equation}
B_{\mu}(0)=\psi_{\mu}^{\prime}(0)=\int_{\mathbb{T}}\zeta\, d\mu(\zeta),\label{eq:2.1}\end{equation}
and that $B_{\delta_{\lambda}}(z)=\lambda$ for all $z\in\mathbb{D}$.
As observed in \cite{SerBerSemi}, \[
\left|B_{\mu}(z)\right|\leq1,\qquad z\in\mathbb{D},\]
and, conversely, any analytic function $B:\,\mathbb{D}\rightarrow\overline{\mathbb{D}}$
is of the form $B_{\mu}$ for a unique probability measure $\mu$
on $\mathbb{T}$. 

Let $\mu_{1}$ and $\mu_{2}$ be two probability measures on $\mathbb{T}$.
As shown in \cite{Franz,BerBoolean}, the multiplicative boolean convolution
$\mu_{1}\utimes\mu_{2}$ is characterized by the following identity\begin{equation}
B_{\mu_{1}\Utimes\mu_{2}}(z)=B_{\mu_{1}}(z)B_{\mu_{2}}(z),\qquad z\in\mathbb{D}.\label{eq:2.2}\end{equation}

It is easy to verify that weak convergence of probability measures
can be translated in terms of the corresponding functions $B$ . More
precisely, given probability measures $\mu$ and $\{\mu_{n}\}_{n=1}^{\infty}$
on $\mathbb{T}$, the sequence $\mu_{n}$ converges weakly to $\mu$
if and only if the sequence $B_{\mu_{n}}(z)$ converges to $B_{\mu}(z)$
uniformly on the compact subsets of $\mathbb{D}$.

A probability measure $\nu$ on $\mathbb{T}$ is \emph{$\Utime$-infinitely
divisible} if, for each $n\in\mathbb{N}$, there exists a probability
measure $\nu_{n}$ on $\mathbb{T}$ such that\[
\nu=\underbrace{\nu_{n}\utimes\nu_{n}\utimes\cdots\utimes\nu_{n}}_{n\:\text{times}}.\]
The notion of infinite divisibility related to other convolutions
is defined analogously. 

The $\Utime$-infinite divisibility is characterized in \cite{Franz}
as follows. A probability measure $\nu$ is $\Utime$-infinitely divisible
if and only if either $\nu$ is Haar measure $m$ (i.e., normalized
arclength measure on $\mathbb{T}$), or the function $B_{\nu}$ can
be expressed as\begin{equation}
B_{\nu}(z)=\gamma\exp\left(-\int_{\mathbb{T}}\frac{1+\zeta z}{1-\zeta z}\, d\sigma(\zeta)\right),\qquad z\in\mathbb{D},\label{eq:2.3}\end{equation}
where $\gamma\in\mathbb{T}$, and $\sigma$ is a finite positive Borel
measure on $\mathbb{T}$. In other words, a measure $\nu$ is $\Utime$-infinitely
divisible if and only if either $B_{\nu}(z)=0$ for all $z\in\mathbb{D}$,
or $0\notin B_{\nu}(\mathbb{D})$. We use the notation $\nu_{\Utimes}^{\gamma,\sigma}$
to denote the $\Utime$-infinitely divisible measure $\nu$ determined
by $\gamma$ and $\sigma$. 

Free multiplicative convolution $\boxtimes$ for probability measures
on the unit circle was introduced by Voiculescu \cite{Voimul}. For
the definition of $\boxtimes$, we refer to \cite{VKN}. Throughout
this paper, we will use the notation $\mathcal{M}_{\mathbb{T}}^{\times}$
to denote the collection of all Borel probability measures $\nu$
on $\mathbb{T}$ with nonzero first moment, i.e., $\int_{\mathbb{T}}\zeta\, d\nu(\zeta)\neq0$. 

In this paper we will require the following characterization \cite{BerVoiBdd}
of $\boxtimes$-infinite divisibility. If a measure $\nu$ is in the
class $\mathcal{M}_{\mathbb{T}}^{\times}$, then the function $\psi_{\nu}$
will have an inverse $\psi^{-1}$ in a neighborhood of zero. In this
case one defines \[
\Sigma_{\nu}(z)=\frac{1}{z}\psi_{\nu}^{-1}\left(\frac{z}{1-z}\right)\]
for $z$ near the origin, and the remarkable identity $\Sigma_{\mu\boxtimes\nu}(z)=\Sigma_{\mu}(z)\Sigma_{\nu}(z)$
holds for $z$ in a neighborhood of zero where three involved functions
are defined. A measure $\nu\in\mathcal{M}_{\mathbb{T}}^{\times}$
is $\boxtimes$-infinitely divisible if and only if the function $\Sigma_{\nu}$
can be expressed as \[
\Sigma_{\nu}(z)=\gamma\exp\left(\int_{\mathbb{T}}\frac{1+\zeta z}{1-\zeta z}\, d\sigma(\zeta)\right),\qquad z\in\mathbb{D},\]
where $\left|\gamma\right|=1$, and $\sigma$ is a finite positive
Borel measure on $\mathbb{T}$. We will use the notation $\nu_{\boxtimes}^{\gamma,\sigma}$
to denote the $\boxtimes$-infinitely divisible measure $\nu$ in
this case. The Haar measure $m$ is the only $\boxtimes$-infinitely
divisible probability measure on $\mathbb{T}$ with zero first moment.

\subsection{Additive boolean and free convolutions on the real line }

Set $\mathbb{C}^{+}=\{ z\in\mathbb{C}:\,\Im z>0\}$ and $\mathbb{C}^{-}=-\mathbb{C}^{+}$.
For $\alpha,\beta>0$, define the cone $\Gamma_{\alpha}=\{ z=x+iy\in\mathbb{C}^{+}:\,\left|x\right|<\alpha y\}$
and the truncated cone $\Gamma_{\alpha,\beta}=\{ z=x+iy\in\Gamma_{\alpha}:\, y>\beta\}$.
We associate every measure $\mu\in\mathcal{M}_{\mathbb{R}}$ its \emph{Cauchy
transform} \[
G_{\mu}(z)=\int_{-\infty}^{\infty}\frac{1}{z-t}\, d\mu(t),\qquad z\in\mathbb{C}^{+},\]
and its reciprocal $F_{\mu}=1/G_{\mu}:\,\mathbb{C}^{+}\rightarrow\mathbb{C}^{+}$.
Then we have $\Im z\leq\Im F_{\mu}(z)$ so that the function $E_{\mu}(z)=z-F_{\mu}(z)$
takes values in $\mathbb{C}^{-}\cup\mathbb{R}$. The function $E_{\mu}$
is such that $E_{\mu}(z)/z\rightarrow0$ as $z\rightarrow\infty$
nontangentially (i.e., $\left|z\right|\rightarrow\infty$ but z stays
within a cone $\Gamma_{\alpha}$ for some $\alpha>0$.) Conversely,
any analytic function $E:\,\mathbb{C}^{+}\rightarrow\mathbb{C}^{-}\cup\mathbb{R}$
such that $E_{\mu}(z)/z\rightarrow0$ as $z\rightarrow\infty$ nontangentially
is of the form $E_{\mu}$ for a unique probability measure $\mu$
on $\mathbb{R}$. 

For $\mu_{1},\mu_{2}\in\mathcal{M}_{\mathbb{R}}$, the additive boolean
convolution $\mu_{1}\uplus\mu_{2}$ is characterized \cite{SpeWoro,BerBoolean}
by the identity\[
E_{\mu_{1}\uplus\mu_{2}}(z)=E_{\mu_{1}}(z)+E_{\mu_{2}}(z),\qquad z\in\mathbb{C}^{+}.\]

Let $\{\mu_{n}\}_{n=1}^{\infty}$ be a sequence in $\mathcal{M}_{\mathbb{R}}$.
As shown in \cite[Proposition 6.2]{BerPataStable}, the sequence $\mu_{n}$
converges weakly to a probability measure $\mu\in\mathcal{M}_{\mathbb{R}}$
if and only if there exists $\beta>0$ such that $\lim_{n\rightarrow\infty}E_{\mu_{n}}(iy)=E_{\mu}(iy)$
for every $y>\beta$, and $E_{\mu_{n}}(iy)=o(y)$ uniformly in $n$
as $y\rightarrow\infty$. 

Every measure $\nu\in\mathcal{M}_{\mathbb{R}}$ is $\uplus$-infinitely
divisible \cite{SpeWoro}. The function $E_{\nu}$ has a Nevanlinna
representation \cite{Achieser}\[
E_{\nu}(z)=\gamma+\int_{-\infty}^{\infty}\frac{1+tz}{z-t}\, d\sigma(t),\qquad z\in\mathbb{C}^{+},\]
where $\gamma\in\mathbb{R}$, and $\sigma$ is a finite positive Borel
measure on $\mathbb{R}$. We use the notation $\nu_{\uplus}^{\gamma,\sigma}$
to denote the ($\uplus$-infinitely divisible) measure $\nu$. 

The additive free convolution $\boxplus$ was first introduced by
Voiculescu \cite{Voiadd} for compactly supported measures on the
real line (then it was extended by Maassen \cite{Maassen} to measures
with finite variance, and by Bercovici and Voiculescu \cite{BerVoiUnBdd}
to the whole class $\mathcal{M}_{\mathbb{R}}$.) The book \cite{VKN}
also contains a detailed description for the theory related to this
convolution. 

We require a result from \cite{BerVoiUnBdd} regarding characterization
of $\boxplus$-infinite divisibility. We have seen earlier that $E_{\mu}(z)/z\rightarrow0$
as $z\rightarrow\infty$ nontangentially for a measure $\mu\in\mathcal{M}_{\mathbb{R}}$.
It follows that for every $\alpha>0$ there exists $\beta=\beta(\mu,\alpha)>0$
such that the function $F_{\mu}$ has an right inverse $F_{\mu}^{-1}$
defined on $\Gamma_{\alpha,\beta}$. The \emph{Voiculescu transform}
\[
\phi_{\mu}(z)=F_{\mu}^{-1}(z)-z,\qquad z\in\Gamma_{\alpha,\beta},\]
linearizes the free convolution in the sense that the identity $\phi_{\mu\boxplus\nu}(z)=\phi_{\mu}(z)+\phi_{\nu}(z)$
holds for $z$ in a truncated cone where all functions involved are
defined. A measure $\nu\in\mathcal{M}_{\mathbb{R}}$ is $\boxplus$-infinitely
divisible if and only if there exist $\gamma\in\mathbb{R}$ and a
finite positive Borel measure $\sigma$ on $\mathbb{R}$ such that
\[
\phi_{\nu}(z)=\gamma+\int_{-\infty}^{\infty}\frac{1+tz}{z-t}\, d\sigma(t),\qquad z\in\mathbb{C}^{+}.\]
We will denote the above measure $\nu$ by $\nu_{\boxplus}^{\gamma,\sigma}$. 

The L\'{e}vy-Hin\v{c}in formula (see \cite{Billingsley}) characterizes
the $*$-infinitely divisible measures in terms of their Fourier transform
as follows: a measure $\rho\in\mathcal{M}_{\mathbb{R}}$ is $*$-infinitely
divisible if and only if there exist $\gamma\in\mathbb{R}$ and a
finite positive Borel measure $\sigma$ on $\mathbb{R}$ such that
the Fourier transform $\widehat{\rho}$ is given by \[
\widehat{\rho}(t)=\exp\left[i\gamma t+\int_{-\infty}^{\infty}\left(e^{itx}-1-\frac{itx}{1+x^{2}}\right)\,\frac{1+x^{2}}{x^{2}}\, d\sigma(x)\right],\quad t\in\mathbb{R},\]
where $\left(e^{itx}-1-\frac{itx}{1+x^{2}}\right)\,\frac{1+x^{2}}{x^{2}}$
is interpreted as $-t^{2}/2$ for $x=0$. The notation $\nu_{*}^{\gamma,\sigma}$
will be used to denote the $*$-infinitely divisible measure determined
by $\gamma$ and $\sigma$.

We will require the following result which was already noted in a
different form in \cite[Lemma 2.3]{BerJcMul}.

\begin{lem}
Consider a sequence of real numbers $\{ r_{n}\}_{n=1}^{\infty}$ and
triangular arrays $\{ z_{nk}\in\mathbb{C}:n\in\mathbb{N},1\leq k\leq k_{n}\}$,
$\{ w_{nk}\in\mathbb{C}:n\in\mathbb{N},1\leq k\leq k_{n}\}$ and $\{ s_{nk}\in\mathbb{R}:n\in\mathbb{N},1\leq k\leq k_{n}\}$.
Suppose that
\begin{enumerate}
\item all $s_{nk}$'s are nonnegative, and \[
\sup_{n\geq1}\sum_{k=1}^{k_{n}}s_{nk}<+\infty;\]

\item $\Re w_{nk}\leq0$ and $\Re z_{nk}\leq0$, for every $n$ and $k$;
\item $z_{nk}=w_{nk}(1+\varepsilon_{nk})$, where the sequence\[
\varepsilon_{n}=\max_{1\leq k\leq k_{n}}\left|\varepsilon_{nk}\right|\]
converges to zero as $n\rightarrow\infty$;
\item there exists a positive constant $M$ such that \[
\left|\Im w_{nk}\right|\leq M\left|\Re w_{nk}\right|+s_{nk},\qquad n\in\mathbb{N},\:1\leq k\leq k_{n}.\]

\end{enumerate}
Then the sequence $\{\exp(ir_{n}+\sum_{k=1}^{k_{n}}z_{nk})\}_{n=1}^{\infty}$
converges if and only if the sequence $\{\exp(ir_{n}+\sum_{k=1}^{k_{n}}w_{nk})\}_{n=1}^{\infty}$
converges. Moreover, the two sequences have the same limit.

\end{lem}
\begin{proof}
From the assumptions on $\{ z_{nk}\}_{n,k}$ and $\{ w_{nk}\}_{n,k}$,
we deduce that \begin{equation}
\left|\sum_{k=1}^{k_{n}}\left[z_{nk}-w_{nk}\right]\right|\leq(1+M)\varepsilon_{n}\left(-\sum_{k=1}^{k_{n}}\Re w_{nk}\right)+\varepsilon_{n}\sum_{k=1}^{k_{n}}s_{nk},\label{eq:2.4}\end{equation}
and \begin{equation}
(1-\varepsilon_{n}-M\varepsilon_{n})\left(-\sum_{k=1}^{k_{n}}\Re w_{nk}\right)\leq\left(-\sum_{k=1}^{k_{n}}\Re z_{nk}\right)+\varepsilon_{n}\sum_{k=1}^{k_{n}}s_{nk},\label{eq:2.5}\end{equation}
for sufficiently large $n$. Suppose that the sequence $\{\exp(ir_{n}+\sum_{k=1}^{k_{n}}z_{nk})\}_{n=1}^{\infty}$
converges to a complex number $z$. If $z=0$, then we have $\lim_{n\rightarrow\infty}\sum_{k=1}^{k_{n}}\Re z_{nk}=-\infty$.
Hence (2.4) implies that $\lim_{n\rightarrow\infty}\sum_{k=1}^{k_{n}}\Re w_{nk}=-\infty$
so that the sequence $\{\exp(ir_{n}+\sum_{k=1}^{k_{n}}w_{nk})\}_{n=1}^{\infty}$
converges to zero as well. If $z\neq0$, then the sequence $\exp(\sum_{k=1}^{k_{n}}\Re z_{nk})$
converges to $\left|z\right|$ as $n\rightarrow\infty$. In particular,
$\sum_{k=1}^{k_{n}}\Re z_{nk}$ is bounded. By (2.4) and (2.5), we
conclude that $\lim_{n\rightarrow\infty}\exp(\sum_{k=1}^{k_{n}}\Re w_{nk})=\left|z\right|$,
and that \[
\lim_{n\rightarrow\infty}\frac{\exp(i\sum_{k=1}^{k_{n}}\Im w_{nk})}{\exp(i\sum_{k=1}^{k_{n}}\Im z_{nk})}=1.\]
Therefore the sequence $\{\exp(ir_{n}+\sum_{k=1}^{k_{n}}w_{nk})\}_{n=1}^{\infty}$
also converges to $z$. The converse implication is proved in the
same way.
\end{proof}

\section{Multiplicative Boolean Convolution on $\mathbb{T}$}

Fix an infinitesimal array $\{\mu_{nk}:n\in\mathbb{N},1\leq k\leq k_{n}\}$
of probability measures on $\mathbb{T}$. For any neighborhood of
zero $\mathcal{V}\subset\mathbb{D}$, it was proved in \cite[Theorem 2.1]{SerBerMulHincin}
that \begin{equation}
\lim_{n\rightarrow\infty}\psi_{\mu_{nk}}(z)=\frac{z}{1-z}\label{eq:3.1}\end{equation}
holds uniformly in $k$ and $z\in\mathcal{V}$. It follows that, as
$n$ tends to infinity, the sequence $B_{\mu_{nk}}(z)$ converges
to $1$ uniformly in $k$ and $z\in\mathcal{V}$. Thus, (2.1) implies
that each $\mu_{nk}$ has nonzero first moment when $n$ is large.
Hence, for our purposes, we will always assume that each member in
such an array belongs to the class $\mathcal{M}_{\mathbb{T}}^{\times}$.
Another application of (3.1) is that the principal branch of $\log B_{\mu_{nk}}(z)$
is defined in $\mathcal{V}$ for large $n$.

Next, we introduce an auxiliary array $\{\mu_{nk}^{\circ}:n\in\mathbb{N},1\leq k\leq k_{n}\}\subset\mathcal{M}_{\mathbb{T}}^{\times}$
as follows. Fix a constant $\tau\in(0,\pi)$. Define the measures
$\mu_{nk}^{\circ}$ by\[
d\mu_{nk}^{\circ}(\zeta)=d\mu_{nk}(b_{nk}\zeta),\]
where the complex numbers $b_{nk}$ are given by\[
b_{nk}=\exp\left(i\int_{\left|\arg\zeta\right|<\tau}\arg\zeta\, d\mu_{nk}(\zeta)\right).\]
Here $\arg\zeta$ is the principal value of the argument of $\zeta$.
Note that the array $\{\mu_{nk}^{\circ}\}_{n,k}$ is again infinitesimal,
and\begin{equation}
\lim_{n\rightarrow\infty}\max_{1\leq k\leq k_{n}}\left|\arg b_{nk}\right|=0.\label{eq:3.2}\end{equation}
We associate each measure $\mu_{nk}^{\circ}$ the function \[
h_{nk}(z)=-i\int_{\mathbb{T}}\Im\zeta\, d\mu_{nk}^{\circ}(\zeta)+\int_{\mathbb{T}}\frac{1+\zeta z}{1-\zeta z}(1-\Re\zeta)\, d\mu_{nk}^{\circ}(\zeta),\qquad z\in\mathbb{D},\]
and observe that $\Re h_{nk}(z)>0$ for all $z\in\mathbb{D}$ unless
the measure $\mu_{nk}^{\circ}=\delta_{1}$. 

\begin{lem}
If $\varepsilon\in(0,1/4)$, then we have, for sufficiently large
$n$, that \[
1-B_{\mu_{nk}^{\circ}}(z)=h_{nk}\left(\overline{b_{nk}}z\right)(1+v_{nk}(z)),\qquad1\leq k\leq k_{n},\]
where $z$ is in $\mathcal{V}_{\varepsilon}=\{ z\in\mathbb{D}:\,\left|z\right|<\varepsilon\}$.
Moreover, we have \[
\lim_{n\rightarrow\infty}\max_{1\leq k\leq k_{n}}\left|v_{nk}(z)\right|=0\]
uniformly on $\mathcal{V}_{\varepsilon}$. 
\end{lem}
\begin{proof}
Applying (3.1) to the array $\{\mu_{nk}^{\circ}\}_{n,k}$, we deduce,
for large $n$, that \[
\frac{z}{1+z}-\frac{\psi_{\mu_{nk}^{\circ}}\left(\frac{z}{1+z}\right)}{1+\psi_{\mu_{nk}^{\circ}}\left(\frac{z}{1+z}\right)}=\frac{1}{(1+z)^{2}}\left[z-\psi_{\mu_{nk}^{\circ}}\left(\frac{z}{1+z}\right)\right](1+u_{nk}(z)),\]
where \[
\lim_{n\rightarrow\infty}\max_{1\le k\leq k_{n}}\left|u_{nk}(z)\right|=0\]
uniformly on $\{ z:\left|z\right|<1/3\}$. Introducing a change of
variable $z\mapsto z/(1-z)$, we obtain \[
z-\frac{\psi_{\mu_{nk}^{\circ}}(z)}{1+\psi_{\mu_{nk}^{\circ}}(z)}=\left[z\int_{\mathbb{T}}\frac{(1-z)(1-\zeta)}{1-\zeta z}\, d\mu_{nk}^{\circ}(\zeta)\right]\left(1+u_{nk}\left(\frac{z}{1-z}\right)\right),\qquad z\in\mathcal{V}_{\varepsilon}.\]
Exploiting the identity\[
\frac{(1-z)(1-\zeta)}{1-\zeta z}=-i\Im\zeta+\frac{1+\zeta z}{1-\zeta z}(1-\Re\zeta),\]
we conclude, for sufficiently large $n$, that\begin{align*}
1-B_{\mu_{nk}^{\circ}}(z) & =\frac{1}{z}\left[z-\frac{\psi_{\mu_{nk}^{\circ}}(z)}{1+\psi_{\mu_{nk}^{\circ}}(z)}\right]\\
 & =h_{nk}(z)\left(1+u_{nk}\left(\frac{z}{1-z}\right)\right),\end{align*}
for all $z\in\mathcal{V}_{\varepsilon}$. 

To prove the result, it suffices to show the following claim: for
every $n$ and $k$, we have \[
h_{nk}\left(\overline{b_{nk}}z\right)=h_{nk}(z)(1+w_{nk}(z)),\]
where $\lim_{n\rightarrow\infty}\max_{1\le k\leq k_{n}}\left|w_{nk}(z)\right|=0$
uniformly in $\mathcal{V}_{\varepsilon}$. If the measure $\mu_{nk}^{\circ}=\delta_{1}$,
then we have $h_{nk}(z)=0$ for all $z\in\mathbb{D}$. In this case,
we define the function $w_{nk}$ to be the zero function in $\mathbb{D}$.
If $\mu_{nk}^{\circ}\neq\delta_{1}$, then we define the function
\[
w_{nk}(z)=\frac{h_{nk}\left(\overline{b_{nk}}z\right)}{h_{nk}(z)}-1,\qquad z\in\mathbb{D}.\]
Observe that\begin{align*}
\left|h_{nk}\left(\overline{b_{nk}}z\right)-h_{nk}(z)\right| & =\left|\left(1-\overline{b_{nk}}\right)\int_{\mathbb{T}}\left[\frac{2\zeta z}{\left(1-\zeta z\right)\left(1-\zeta\overline{b_{nk}}z\right)}\right](1-\Re\zeta)\, d\mu_{nk}^{\circ}(\zeta)\right|\\
 & \leq\left|1-\overline{b_{nk}}\right|\int_{\mathbb{T}}\left|\frac{2\zeta z}{\left(1-\zeta z\right)\left(1-\zeta\overline{b_{nk}}z\right)}\right|(1-\Re\zeta)\, d\mu_{nk}^{\circ}(\zeta)\\
 & \leq\frac{2}{(1-\varepsilon)^{2}}\left|1-\overline{b_{nk}}\right|\int_{\mathbb{T}}(1-\Re\zeta)\, d\mu_{nk}^{\circ}(\zeta),\end{align*}
for $z\in\mathcal{V}_{\varepsilon}$ and $\zeta\in\mathbb{T}$. Meanwhile,
Harnack's inequality implies that there exists $L=L(\varepsilon)>0$
such that \[
\left|\Re\left[\frac{1+\zeta z}{1-\zeta z}\right]\right|=\Re\left[\frac{1+\zeta z}{1-\zeta z}\right]\geq L,\qquad z\in\mathcal{V}_{\varepsilon},\:\zeta\in\mathbb{T}.\]
Thus, we have \begin{align*}
\left|h_{nk}(z)\right| & \geq\Re h_{nk}(z)\\
 & =\int_{\mathbb{T}}\Re\left[\frac{1+\zeta z}{1-\zeta z}\right](1-\Re\zeta)\, d\mu_{nk}^{\circ}(\zeta)\\
 & \geq L\int_{\mathbb{T}}(1-\Re\zeta)\, d\mu_{nk}^{\circ}(\zeta).\end{align*}
Combining the above inequalities, we get\begin{align*}
\left|w_{nk}(z)\right| & \leq\frac{\left|h_{nk}\left(\overline{b_{nk}}z\right)-h_{nk}(z)\right|}{\left|h_{nk}(z)\right|}\\
 & \leq\frac{2}{(1-\varepsilon)^{2}L}\left|1-\overline{b_{nk}}\right|\leq\frac{2}{(1-\varepsilon)^{2}L}\left|\arg b_{nk}\right|,\end{align*}
for $z\in\mathcal{V}_{\varepsilon}$. Hence the claim is proved by
(3.2).
\end{proof}
A crucial property for the functions $h_{nk}(z)$ proved in \cite[Lemma 4.1]{BerJcMul}
is that for every neighborhood of zero $\mathcal{V}\subset\mathbb{D}$
there exists a constant $M=M(\mathcal{V},\tau)>0$ such that \begin{equation}
\left|\Im h_{nk}(z)\right|\leq M\Re h_{nk}(z),\qquad z\in\mathcal{V},\;1\leq k\leq k_{n},\label{eq:3.3}\end{equation}
for sufficiently large $n$.

\begin{prop}
Let $\{\lambda_{n}\}_{n=1}^{\infty}$be a sequence in $\mathbb{T}$.
Suppose that $B:\,\mathbb{D}\rightarrow\overline{\mathbb{D}}$ is
an analytic function. Then \[
\lim_{n\rightarrow\infty}\lambda_{n}\prod_{k=1}^{k_{n}}B_{\mu_{nk}}(z)=B(z)\]
uniformly on the compact subsets of $\mathbb{D}$ if and only if \[
\lim_{n\rightarrow\infty}\exp\left(i\arg\lambda_{n}+i\sum_{k=1}^{k_{n}}\arg b_{nk}-\sum_{k=1}^{k_{n}}h_{nk}(z)\right)=B(z)\]
uniformly on the compact subsets of $\mathbb{D}$. 
\end{prop}
\begin{proof}
Suppose that the sequence $\{\exp(i\arg\lambda_{n}+i\sum_{k=1}^{k_{n}}\arg b_{nk}-\sum_{k=1}^{k_{n}}h_{nk}(z))\}_{n=1}^{\infty}$
converges to $B(z)$ uniformly on the compact subsets of $\mathbb{D}$.
Note that\[
\log w=w-1+o(\left|w-1\right|)\]
as $w\rightarrow1$, and $B_{\mu_{nk}^{\circ}}(z)=\overline{b_{nk}}B_{\mu_{nk}}\left(\overline{b_{nk}}z\right)$
for every $z\in\mathbb{D}$. For $z$ near the origin, Lemma 3.1 shows
that $\log\left(\overline{b_{nk}}B_{\mu_{nk}}(z)\right)=-h_{nk}(z)(1+o(1))$
uniformly in $k$ as $n$ tends to infinity. Then Lemma 2.1 implies
that \[
\lim_{n\rightarrow\infty}\lambda_{n}\prod_{k=1}^{k_{n}}B_{\mu_{nk}}(z)=B(z)\]
uniformly in a neighborhood of zero. Moreover, this convergence is
actually uniform on the compact subsets of $\mathbb{D}$ since the
family $\{\lambda_{n}\prod_{k=1}^{k_{n}}B_{\mu_{nk}}(z)\}_{n=1}^{\infty}$
is normal. The converse implication is proved in the same way. 
\end{proof}
\begin{lem}
Let $\{\nu_{n}\}_{n=1}^{\infty}$ be a sequence of $\Utime$-infinitely
divisible measures on $\mathbb{T}$. If the sequence $\nu_{n}$ converges
weakly to a probability measure $\nu$, then the measure $\nu$ is
$\Utime$-infinitely divisible.
\end{lem}
\begin{proof}
The weak convergence of $\nu_{n}$ implies that the sequence $B_{\nu_{n}}(z)$
converges to $B_{\nu}(z)$ uniformly on the compact subsets of $\mathbb{D}$.
If the function $B_{\nu}$ is nonvanishing in $\mathbb{D}$, then
the measure $\nu$ is $\Utime$-infinitely divisible. On the other
hand, if $B_{\nu}(z_{0})=0$ for some $z_{0}\in\mathbb{D}$, then
Rouch\'{e}'s theorem implies that there exists an $N=N(z_{0})\in\mathbb{N}$
such that the function $B_{\nu_{n}}(z)$ also has a zero in the disk
$\{ z:\left|z-z_{0}\right|<1-\left|z_{0}\right|\}$ whenever $n\geq N$.
Since each $\nu_{n}$ is $\Utime$-infinitely divisible, we conclude
in this case that $\nu_{n}$ is the Haar measure $m$ for all $n\geq N$.
Consequently, the measure $\nu$ must be $m$ as well. 
\end{proof}
Our next result is the boolean analogue of Hin\v{c}in's theorem. 

\begin{thm}
Let $\{\lambda_{n}\}_{n=1}^{\infty}$ be a sequence in $\mathbb{T}$.
If the sequence of measures \[
\delta_{\lambda_{n}}\utimes\mu_{n1}\utimes\mu_{n2}\utimes\cdots\utimes\mu_{nk_{n}}\]
 converges weakly on $\mathbb{T}$ to a probability measure $\nu$,
then $\nu$ is $\Utime$-infinitely divisible. 
\end{thm}
\begin{proof}
From (2.2) and the weak convergence of $\delta_{\lambda_{n}}\utimes\mu_{n1}\utimes\mu_{n2}\utimes\cdots\utimes\mu_{nk_{n}}$,
we have \[
\lim_{n\rightarrow\infty}\lambda_{n}\prod_{k=1}^{k_{n}}B_{\mu_{nk}}(z)=B_{\nu}(z)\]
uniformly on the compact subsets of $\mathbb{D}$. Observe that the
function $-\sum_{k=1}^{k_{n}}h_{nk}(z)$ has negative real part in
$\mathbb{D}$, and hence there exists a $\Utime$-infinitely divisible
measure $\nu_{n}$ on $\mathbb{T}$ such that \[
B_{\nu_{n}}(z)=\exp\left(i\arg\lambda_{n}+i\sum_{k=1}^{k_{n}}\arg b_{nk}-\sum_{k=1}^{k_{n}}h_{nk}(z)\right),\qquad z\in\mathbb{D}.\]
Proposition 3.2 then implies that the sequence $\nu_{n}$ converges
weakly to $\nu$. The $\Utime$-infinitely divisibility of the measure
$\nu$ follows immediately by Lemma 3.3.
\end{proof}
Fix $\gamma\in\mathbb{T}$ and a finite positive Borel measure $\sigma$
on $\mathbb{T}$. 

\begin{thm}
For the infinitesimal array $\{\mu_{nk}\}_{n,k}\subset\mathcal{M}_{\mathbb{T}}^{\times}$
and a sequence $\{\lambda_{n}\}_{n=1}^{\infty}\subset\mathbb{T}$,
the following statements are equivalent:
\begin{enumerate}
\item The sequence $\delta_{\lambda_{n}}\utimes\mu_{n1}\utimes\mu_{n2}\utimes\cdots\utimes\mu_{nk_{n}}$
converges weakly to $\nu_{\Utimes}^{\gamma,\sigma}$.
\item The sequence $\delta_{\lambda_{n}}\boxtimes\mu_{n1}\boxtimes\mu_{n2}\boxtimes\cdots\boxtimes\mu_{nk_{n}}$
converges weakly to $\nu_{\boxtimes}^{\gamma,\sigma}$.
\item The sequence of measures \[
d\sigma_{n}(\zeta)=\sum_{k=1}^{k_{n}}(1-\Re\zeta)\, d\mu_{nk}^{\circ}(\zeta)\]
converges weakly on $\mathbb{T}$ to $\sigma$, and the limit \[
\lim_{n\rightarrow\infty}\gamma_{n}=\gamma\]
exists, where\[
\gamma_{n}=\exp\left(i\arg\lambda_{n}+i\sum_{k=1}^{k_{n}}\arg b_{nk}+i\sum_{k=1}^{k_{n}}\int_{\mathbb{T}}\Im\zeta\, d\mu_{nk}^{\circ}(\zeta)\right).\]

\end{enumerate}
\end{thm}
\begin{proof}
The equivalence of (2) and (3) has been proved in \cite{BerJcMul}.
We will focus on the equivalence of (1) and (3). Assume that (1) holds.
Then we have \[
\lim_{n\rightarrow\infty}\lambda_{n}\prod_{k=1}^{k_{n}}B_{\mu_{nk}}(z)=B_{\nu_{\Utimes}^{\gamma,\sigma}}(z)=\gamma\exp\left(-\int_{\mathbb{T}}\frac{1+\zeta z}{1-\zeta z}\, d\sigma(\zeta)\right)\]
uniformly on the compact subsets of $\mathbb{D}$. Proposition 3.2
then shows that\begin{equation}
\lim_{n\rightarrow\infty}\exp\left(i\arg\lambda_{n}+i\sum_{k=1}^{k_{n}}\arg b_{nk}-\sum_{k=1}^{k_{n}}h_{nk}(z)\right)=\gamma\exp\left(-\int_{\mathbb{T}}\frac{1+\zeta z}{1-\zeta z}\, d\sigma(\zeta)\right)\label{eq:3.4}\end{equation}
uniformly on the compact subsets of $\mathbb{D}$. Taking the absolute
value on both sides, we conclude that \begin{equation}
\lim_{n\rightarrow\infty}\exp\left(-\sum_{k=1}^{k_{n}}\Re h_{nk}(z)\right)=\exp\left(-\int_{\mathbb{T}}\Re\left[\frac{1+\zeta z}{1-\zeta z}\right]\, d\sigma(\zeta)\right),\qquad z\in\mathbb{D}.\label{eq:3.5}\end{equation}
Since \[
\exp\left(i\arg\lambda_{n}+i\sum_{k=1}^{k_{n}}\arg b_{nk}-\sum_{k=1}^{k_{n}}h_{nk}(z)\right)=\gamma_{n}\exp\left(-\int_{\mathbb{T}}\frac{1+\zeta z}{1-\zeta z}\, d\sigma_{n}(\zeta)\right),\]
and the real part of the function $\sum_{k=1}^{k_{n}}h_{nk}(z)$ is
the Poisson integral of the measure $d\sigma_{n}\left(\overline{\zeta}\right)$,
the equation (3.5) uniquely determines the measure $\sigma$ which
is the weak cluster point of $\{\sigma_{n}\}_{n=1}^{\infty}$. Hence,
$\sigma_{n}$ must converge weakly to $\sigma$. The convergence property
of the sequence $\gamma_{n}$ follows immediately by letting $z=0$
in (3.4) and (3.5).

For the converse implication from (3) to (1), one can easily reverse
the above steps to reach (1) by Proposition 3.2. The details are left
to the reader. 
\end{proof}
The equivalent condition for the weak convergence of $\delta_{\lambda_{n}}\boxtimes\mu_{n1}\boxtimes\mu_{n2}\boxtimes\cdots\boxtimes\mu_{nk_{n}}$
to Haar measure $m$ was given in \cite[Theorem 4.4]{BerJcMul}. It
turns out that the same condition is also equivalent to the weak convergence
of $\delta_{\lambda_{n}}\utimes\mu_{n1}\utimes\mu_{n2}\utimes\cdots\utimes\mu_{nk_{n}}$
to $m$. We will not provide the details of this proof because they
are entirely analogous to the free case. We only point out the relevant
fact needed in the proof is that\[
B_{\mu}(0)=\int_{\mathbb{T}}\zeta\, d\mu(\zeta)\]
for all probability measure $\mu$ on $\mathbb{T}$. 

\begin{thm}
For the infinitesimal array $\{\mu_{nk}\}_{n,k}\subset\mathcal{M}_{\mathbb{T}}^{\times}$
and a sequence $\{\lambda_{n}\}_{n=1}^{\infty}\subset\mathbb{T}$,
the following statements are equivalent:
\begin{enumerate}
\item The sequence $\delta_{\lambda_{n}}\utimes\mu_{n1}\utimes\mu_{n2}\utimes\cdots\utimes\mu_{nk_{n}}$
converges weakly to $m$.
\item The sequence $\delta_{\lambda_{n}}\boxtimes\mu_{n1}\boxtimes\mu_{n2}\boxtimes\cdots\boxtimes\mu_{nk_{n}}$
converges weakly to $m$.
\item \[
\lim_{n\rightarrow\infty}\sum_{k=1}^{k_{n}}\int_{\mathbb{T}}(1-\Re\zeta)\, d\mu_{nk}^{\circ}(\zeta)=+\infty.\]

\end{enumerate}
\end{thm}
We conclude this section by using Theorem 3.5 to determine the multiplicative
boolean analogues of Gaussian and Poisson laws on $\mathbb{R}$. The
following result generates a measure analogous to the Gaussian distribution
on the real line. 

\begin{cor}
For every $t>0$, the function\[
B(z)=\exp\left(-\frac{t}{2}\left(\frac{1+z}{1-z}\right)\right),\qquad z\in\mathbb{D},\]
is of the form $B=B_{\nu}$ for some $\Utime$-infinitely divisible
measure $\nu\in\mathcal{M}_{\mathbb{T}}^{\times}$.
\end{cor}
\begin{proof}
For $n>t$, we define \[
\mu_{nk}=\mu_{n}=\frac{1}{2}\left(\delta_{\xi_{n}}+\delta_{\overline{\xi_{n}}}\right),\qquad1\leq k\leq n,\]
where \[
\xi_{n}=\sqrt{1-\frac{t}{n}}+i\sqrt{\frac{t}{n}}.\]
To apply Theorem 3.5, we choose $\tau=1$ so that $b_{nk}=1$ for
every $n$ and $k$. Hence we have $\mu_{n}^{\circ}=\mu_{n}$. As
in the statement of Theorem 3.5, we define the measures\[
d\sigma_{n}(\zeta)=n(1-\Re\zeta)\, d\mu_{n}(\zeta),\]
and the numbers $\gamma_{n}=\exp\left(in\int_{\mathbb{T}}\Im\zeta\, d\mu_{n}(\zeta)\right)$.
Note that $\gamma_{n}=1$ for all $n\in\mathbb{N}$, and the $p$-th
Fourier coefficient $\widehat{\sigma_{n}}(p)$ of the measure $\sigma_{n}$
is given by \[
\widehat{\sigma_{n}}(p)=\int_{\mathbb{T}}\zeta^{^{p}}n(1-\Re\zeta)\, d\mu_{n}(\zeta)=n\Re\xi_{n}^{^{p}}(1-\Re\xi_{n}),\]
where $p$ is an integer. Since $\lim_{n\rightarrow\infty}\widehat{\sigma_{n}}(p)=t/2$
for all $p$, we conclude that the sequence $\sigma_{n}$ converges
weakly on $\mathbb{T}$ to the measure \[
\sigma=\frac{t}{2}\delta_{1}.\]
Theorem 3.5 then implies that the sequence $\underbrace{\mu_{n}\utimes\mu_{n}\utimes\cdots\utimes\mu_{n}}_{n\:\text{times}}$
converges weakly to $\nu_{\Utimes}^{1,\sigma}$ as $n\rightarrow\infty$.
The desired result now follows from (2.3).
\end{proof}
\begin{defn*}
A $\Utime$-infinitely divisible measure $\nu_{\Utimes}^{\gamma,\sigma}\in\mathcal{M}_{\mathbb{T}}^{\times}$
is said to be \emph{$\Utime$-normal} if the measure $\sigma$ is
concentrated in the point $1$ (i.e., $\sigma=\sigma(\mathbb{T})\delta_{1}$).
\end{defn*}
Our next result produces a boolean analogue of the Poisson distribution
on $\mathbb{R}$.

\begin{cor}
For every $t>0$ and $\lambda\in\mathbb{T}$, the function \[
B(z)=\exp\left(-t(1-\lambda)\left(\frac{1-z}{1-\lambda z}\right)\right),\qquad z\in\mathbb{D},\]
is of the form $B=B_{\nu}$ for some $\Utime$-infinitely divisible
measure $\nu\in\mathcal{M}_{\mathbb{T}}^{\times}$.
\end{cor}
\begin{proof}
Note that $B=B_{\delta_{1}}$ when $\lambda=1$. Assume now $\lambda\neq1$.
This time we set \[
\mu_{nk}=\mu_{n}=\left(1-\frac{t}{n}\right)\delta_{1}+\frac{t}{n}\delta_{\lambda},\qquad1\leq k\leq n,\]
and we choose $\tau=\left|\arg\lambda\right|/2$ so that $\mu_{n}^{\circ}=\mu_{n}$.
Meanwhile, we define the measures $\sigma_{n}$ and the numbers $\gamma_{n}$
as in the proof of Corollary 3.7. Then we have $\widehat{\sigma_{n}}(p)=t\lambda^{^{p}}(1-\Re\lambda)$
and $\gamma_{n}=e^{it\Im\lambda}$, for all $p\in\mathbb{Z}$ and
$n\in\mathbb{N}$. Thus, the measures $\sigma_{n}$ converge weakly
on $\mathbb{T}$ to the measure \[
\sigma=t(1-\Re\lambda)\delta_{\lambda},\]
while the number $\gamma=e^{it\Im\lambda}$. Then the proof is completed
by Theorem 3.5 and the following observation:\begin{align*}
it\Im\lambda-t(1-\Re\lambda)\frac{1+\lambda z}{1-\lambda z} & =-t\left[-i\Im\lambda+(1-\Re\lambda)\frac{1+\lambda z}{1-\lambda z}\right]\\
 & =-t\left[\frac{(1-\lambda)(1-z)}{1-\lambda z}\right].\end{align*}

\end{proof}

\section{Classical Convolution on $\mathbb{T}$ }

Consider an infinitesimal array $\{\mu_{nk}\}_{n,k}\subset\mathcal{M}_{\mathbb{T}}^{\times}$
and a sequence $\{\lambda_{n}\}_{n=1}^{\infty}\subset\mathbb{T}$,
we define \[
\mu_{n}=\delta_{\lambda_{n}}\utimes\mu_{n1}\utimes\mu_{n2}\utimes\cdots\utimes\mu_{nk_{n}},\]
and\[
\nu_{n}=\delta_{\lambda_{n}}\circledast\mu_{n1}\circledast\mu_{n2}\circledast\cdots\circledast\mu_{nk_{n}},\]
for every $n\in\mathbb{N}$. The aim of current section is to investigate
connections between the asymptotic distributional behavior of $\{\mu_{n}\}_{n=1}^{\infty}$
and that of $\{\nu_{n}\}_{n=1}^{\infty}$. For our purposes, we introduce
the complex numbers\[
b_{nk}=\exp\left(i\int_{\left|\arg\zeta\right|<1}\arg\zeta\, d\mu_{nk}(\zeta)\right),\]
and the centered measures $d\mu_{nk}^{\circ}(\zeta)=d\mu_{nk}(b_{nk}\zeta)$.
Note that we have \[
\widehat{\mu_{nk}}(p)=b_{nk}^{p}\widehat{\mu_{nk}^{\circ}}(p)\]
for any integer $p$, and that the function \[
\frac{\zeta^{^{p}}-1-ip\Im\zeta}{1-\Re\zeta}\]
is continuous and bounded on $\mathbb{T}$. (The value of this function
for $\zeta=1$ is set at $-p^{2}$ in order to preserve its continuity
at that point.) 

\begin{thm}
Assume that $\gamma\in\mathbb{T}$, and that $\sigma$ is a finite
positive Borel measure on $\mathbb{T}$. If the sequence $\mu_{n}$
converges weakly to $\nu_{\Utimes}^{\gamma,\sigma}$ , then there
exists a probability measure $\nu$ on $\mathbb{T}$ such that the
sequence $\nu_{n}$ converges weakly to $\nu$. Moreover, the Fourier
coefficients of the limit law $\nu$ can be calculated by the formula:
\begin{equation}
\widehat{\nu}(p)=\gamma^{p}\exp\left(\int_{\mathbb{T}}\frac{\zeta^{^{p}}-1-ip\Im\zeta}{1-\Re\zeta}\, d\sigma(\zeta)\right),\qquad p\in\mathbb{Z}.\label{eq:4.1}\end{equation}

\end{thm}
\begin{proof}
Observe that $\widehat{\nu_{n}}(0)=1$ for all $n\in\mathbb{N}$,
and the right side of (4.1) is $1$ when $p=0$. Fix now a nonzero
integer $p$. To prove the theorem, it suffices to show that the sequence
$\{\widehat{\nu_{n}}(p)\}_{n=1}^{\infty}$ has a limit, and that this
limit can be identified as the right side of (4.1). Since the array
$\{\mu_{nk}^{\circ}\}_{n,k}$ is infinitesimal, the principal logarithm
of $\widehat{\mu_{nk}^{\circ}}(p)$ exists when $n$ is sufficiently
large. Moreover, we have \begin{equation}
\widehat{\nu_{n}}(p)=\exp\left(ip\arg\lambda_{n}+ip\sum_{k=1}^{k_{n}}\arg b_{nk}+\sum_{k=1}^{k_{n}}\log\widehat{\mu_{nk}^{\circ}}(p)\right)\label{eq:4.2}\end{equation}
for large $n$. Define the complex numbers $A_{nk}=A_{nk}(p)=\widehat{\mu_{nk}^{\circ}}(p)-1$,
and set \[
d\sigma_{n}(\zeta)=\sum_{k=1}^{k_{n}}(1-\Re\zeta)\, d\mu_{nk}^{\circ}(\zeta),\]
and \[
\gamma_{n}=\exp\left(i\arg\lambda_{n}+i\sum_{k=1}^{k_{n}}\arg b_{nk}+i\sum_{k=1}^{k_{n}}\int_{\mathbb{T}}\Im\zeta\, d\mu_{nk}^{\circ}(\zeta)\right).\]
By Theorem 3.5, the measures $\sigma_{n}$ converge weakly on $\mathbb{T}$
to the measure $\sigma$, and the limit of the sequence $\gamma_{n}$
is $\gamma$. Note that\begin{align*}
\exp\left(ip\arg\lambda_{n}+ip\sum_{k=1}^{k_{n}}\arg b_{nk}+\sum_{k=1}^{k_{n}}A_{nk}\right) & =\gamma_{n}^{p}\exp\left(\sum_{k=1}^{k_{n}}\left[A_{nk}-\int_{\mathbb{T}}ip\Im\zeta\, d\mu_{nk}^{\circ}(\zeta)\right]\right)\\
 & =\gamma_{n}^{p}\exp\left(\sum_{k=1}^{k_{n}}\int_{\mathbb{T}}\zeta^{^{p}}-1-ip\Im\zeta\, d\mu_{nk}^{\circ}(\zeta)\right)\\
 & =\gamma_{n}^{p}\exp\left(\int_{\mathbb{T}}\frac{\zeta^{^{p}}-1-ip\Im\zeta}{1-\Re\zeta}\, d\sigma_{n}(\zeta)\right).\end{align*}
Therefore, we deduce that \[
\lim_{n\rightarrow\infty}\exp\left(ip\arg\lambda_{n}+ip\sum_{k=1}^{k_{n}}\arg b_{nk}+\sum_{k=1}^{k_{n}}A_{nk}\right)=\gamma^{p}\exp\left(\int_{\mathbb{T}}\frac{\zeta^{^{p}}-1-ip\Im\zeta}{1-\Re\zeta}\, d\sigma(\zeta)\right).\]

The infinitesimality of the array $\{\mu_{nk}^{\circ}\}_{n,k}$ implies
that $\max_{1\leq k\leq k_{n}}\left|A_{nk}\right|\rightarrow0$ as
$n\rightarrow\infty$. Hence for sufficiently large $n$ the expansion\[
\log\widehat{\mu_{nk}^{\circ}}(p)=\log(1+A_{nk})=A_{nk}-\frac{1}{2}A_{nk}^{2}+\frac{1}{3}A_{nk}^{3}-\cdots\]
holds. Thus, we deduce that $\log\widehat{\mu_{nk}^{\circ}}(p)=A_{nk}(1+o(1))$
uniformly in $k$ as $n\rightarrow\infty$. 

Denote by $\mathcal{U}_{p}$ the set of all complex numbers $\zeta\in\mathbb{T}$
such that $3\left|\arg\zeta\right|<\min\{1,\left|\pi/p\right|\}$,
and by $\mathcal{V}_{p}$ the set of all $\zeta\in\mathcal{U}_{p}$
such that $6\left|\arg\zeta\right|<\min\{1,\left|\pi/p\right|\}$.
We also introduce the sets $\mathcal{U}_{p}^{\circ}=\mathcal{U}_{p}^{\circ}(n,k)=\{ b_{nk}\zeta:\,\zeta\in\mathcal{U}_{p}\}$.
By (3.2), we have \begin{align*}
\left|\int_{\mathcal{U}_{p}}\arg\zeta\, d\mu_{nk}^{\circ}(\zeta)\right| & =\left|\arg b_{nk}-\int_{\{\left|\arg\zeta\right|<1\}\setminus\mathcal{U}_{p}^{\circ}}\arg\zeta\, d\mu_{nk}(\zeta)-\arg b_{nk}\int_{\mathcal{U}_{p}^{\circ}}\, d\mu_{nk}(\zeta)\right|\\
 & =\left|\arg b_{nk}\mu_{nk}^{\circ}(\mathbb{T}\setminus\mathcal{U}_{p})-\int_{\{\left|\arg\zeta\right|<1\}\setminus\mathcal{U}_{p}^{\circ}}\arg\zeta\, d\mu_{nk}(\zeta)\right|\\
 & \leq2\mu_{nk}^{\circ}(\mathbb{T}\setminus\mathcal{V}_{p}),\end{align*}
for sufficiently large $n$. Hence we conclude, for large $n$, that\begin{align*}
\left|\Im A_{nk}\right| & \leq\int_{\mathcal{U}_{p}}\left|\Im\zeta^{^{p}}-p\arg\zeta\right|\, d\mu_{nk}^{\circ}(\zeta)+\left|\int_{\mathcal{U}_{p}}p\arg\zeta\, d\mu_{nk}^{\circ}(\zeta)\right|+\int_{\mathbb{T}\setminus\mathcal{U}_{p}}\left|\Im\zeta^{^{p}}\right|\, d\mu_{nk}^{\circ}(\zeta)\\
 & \leq2\int_{\mathcal{U}_{p}}(1-\Re\zeta^{^{p}})\, d\mu_{nk}^{\circ}(\zeta)+(2\left|p\right|+1)\mu_{nk}^{\circ}(\mathbb{T}\setminus\mathcal{V}_{p})\\
 & \leq2\left|\Re A_{nk}\right|+(2\left|p\right|+1)\mu_{nk}^{\circ}(\mathbb{T}\setminus\mathcal{V}_{p}).\end{align*}
Meanwhile, the weak convergence of $\sigma_{n}$ implies that \[
\lim_{n\rightarrow\infty}\int_{\mathbb{T}\setminus\mathcal{V}_{p}}\frac{1}{1-\Re\zeta}\, d\sigma_{n}(\zeta)=\int_{\mathbb{T}\setminus\mathcal{V}_{p}}\frac{1}{1-\Re\zeta}\, d\sigma(\zeta).\]
Since \[
\sum_{k=1}^{k_{n}}\mu_{nk}^{\circ}(\mathbb{T}\setminus\mathcal{V}_{p})=\int_{\mathbb{T}\setminus\mathcal{V}_{p}}\frac{1}{1-\Re\zeta}\, d\sigma_{n}(\zeta),\]
we conclude that $\sum_{k=1}^{k_{n}}\mu_{nk}^{\circ}(\mathbb{T}\setminus\mathcal{V}_{p})$
is bounded.

Applying Lemma 2.1 to the arrays $\{ A_{nk}\}_{n,k}$ and $\{\log\widehat{\mu_{nk}^{\circ}}(p)\}_{n,k}$,
we conclude at once that the sequence $\widehat{\nu_{n}}(p)$ converges,
and\begin{align*}
\lim_{n\rightarrow\infty}\widehat{\nu_{n}}(p) & =\lim_{n\rightarrow\infty}\exp\left(ip\arg\lambda_{n}+ip\sum_{k=1}^{k_{n}}\arg b_{nk}+\sum_{k=1}^{k_{n}}A_{nk}\right)\\
 & =\gamma^{p}\exp\left(\int_{\mathbb{T}}\frac{\zeta^{^{p}}-1-ip\Im\zeta}{1-\Re\zeta}\, d\sigma(\zeta)\right).\end{align*}

\end{proof}
\begin{rem*}
Note that (4.1) implies that the limit law $\nu$ in Theorem 4.1 is
$\circledast$-infinitely divisible. Indeed, for every $n\in\mathbb{N}$,
there exists a probability measure $\nu_{n}$ on $\mathbb{T}$ such
that \[
\widehat{\nu_{n}}(p)=\gamma^{\frac{p}{n}}\exp\left(\frac{1}{n}\int_{\mathbb{T}}\frac{\zeta^{^{p}}-1-ip\Im\zeta}{1-\Re\zeta}\, d\sigma(\zeta)\right),\qquad p\in\mathbb{Z}.\]
It follows that $\nu=\underbrace{\nu_{n}\circledast\nu_{n}\circledast\cdots\circledast\nu_{n}}_{n\:\text{times}}$,
and hence the measure $\nu$ is $\circledast$-infinitely divisible. 
\end{rem*}
Suppose $a\in\mathbb{R}$ and $t>0$. Denote by $N(a,t)$ the Gaussian
distribution on $\mathbb{R}$ with mean $a$ and variance $t$, that
is, \[
dN(a,t)(x)=\frac{1}{\sqrt{2\pi t}}e^{-\frac{1}{2t}(x-a)^{2}}\, dx,\qquad-\infty<x<\infty.\]
Let $\tau$ be the continuous homomorphism $x\mapsto e^{ix}$ from
$\mathbb{R}$ into the circle $\mathbb{T}$. A probability measure
$\nu$ on $\mathbb{T}$ is called a normal distribution \cite[Chapter V, Section 5.2]{Heyer}
if $\nu$ is the push-forward measure of a Gaussian law $N(a,t)$
through the map $\tau$. One computes its measure $\nu(S)$ of a Borel
measurable set $S\subset\mathbb{T}$ as\[
\nu(S)=\int_{\arg S}\sum_{n\in\mathbb{Z}}\frac{1}{\sqrt{2\pi t}}e^{-\frac{1}{2t}(u-a+2n\pi)^{2}}\, du,\]
where the set $\arg S=\{\arg\zeta:\,\zeta\in S\}$. Note that $\nu$
is normal if and only if \[
\widehat{\nu}(p)=\exp\left(iap-\frac{t}{2}p^{2}\right),\qquad p\in\mathbb{Z}.\]
It follows that each normal distribution on $\mathbb{T}$ is $\circledast$-infinitely
divisible. The next result shows that the boolean (or free) central
limit theorem holds if and only if the classical central limit theorem
holds. Recall a $\Utime$-normal distribution on $\mathbb{T}$ is
a $\Utime$-infinitely divisible measure $\nu_{\Utimes}^{\gamma,\sigma}$
such that the measure $\sigma$ is concentrated in the point $1$.

\begin{cor}
The sequence $\mu_{n}$ converges weakly on $\mathbb{T}$ to a $\Utime$-normal
distribution if and only if the sequence $\nu_{n}$ converges weakly
on $\mathbb{T}$ to a normal distribution. 
\end{cor}
\begin{proof}
If the sequence $\mu_{n}$ converges weakly to a $\Utime$-normal
distribution $\nu_{\Utimes}^{\gamma,\sigma}$, then Theorem 4.1 shows
that the sequence $\nu_{n}$ converges weakly to a probability measure
$\nu$ such that $\widehat{\nu}(p)=\gamma^{p}\exp(-\sigma(\{1\})p^{2})$
for all $p\in\mathbb{Z}$. Therefore, the measure $\nu$ is a normal
distribution on $\mathbb{T}$. 

Assume now that the sequence $\nu_{n}$ converges weakly on $\mathbb{T}$
to a normal distribution $\nu$. Then we have \[
\lim_{n\rightarrow\infty}\widehat{\nu_{n}}(p)=\widehat{\nu}(p)=\gamma^{p}\exp\left(-\frac{t}{2}p^{2}\right),\qquad p\in\mathbb{Z},\ \]
for some $\gamma\in\mathbb{T}$ and $t>0$. Define the complex numbers
$A_{nk}(p)$, $\gamma_{n}$ and the measures $\sigma_{n}$ as in the
proof of Theorem 4.1, and note that \[
\exp\left(i\arg\lambda_{n}+i\sum_{k=1}^{k_{n}}\arg b_{nk}+\sum_{k=1}^{k_{n}}A_{nk}(1)\right)=\gamma_{n}e^{-\sigma_{n}(\mathbb{T})},\qquad n\in\mathbb{N}.\]
Let us recall, from Section 3, the definition of functions \[
h_{nk}(z)=-i\int_{\mathbb{T}}\Im\zeta\, d\mu_{nk}^{\circ}(\zeta)+\int_{\mathbb{T}}\frac{1+\zeta z}{1-\zeta z}(1-\Re\zeta)\, d\mu_{nk}^{\circ}(\zeta),\qquad z\in\mathbb{D},\]
and observe that $\left|\Im A_{nk}(1)\right|=\left|\Im h_{nk}(0)\right|$
and $\left|\Re A_{nk}(1)\right|=\left|\Re h_{nk}(0)\right|$. Then
(3.3) shows that there exists $M>0$ such that $\left|\Im A_{nk}(1)\right|\leq M\left|\Re A_{nk}(1)\right|$
for large $n$. Since $\log\widehat{\mu_{nk}^{\circ}}(1)=A_{nk}(1)(1+o(1))$
uniformly in $k$ as $n\rightarrow\infty$, Lemma 2.1 and (4.2) imply
that \begin{align*}
\lim_{n\rightarrow\infty}\gamma_{n}e^{-\sigma_{n}(\mathbb{T})} & =\lim_{n\rightarrow\infty}\exp\left(i\arg\lambda_{n}+i\sum_{k=1}^{k_{n}}\arg b_{nk}+\sum_{k=1}^{k_{n}}A_{nk}(1)\right)\\
 & =\lim_{n\rightarrow\infty}\widehat{\nu_{n}}(1)=\gamma\exp\left(-\frac{t}{2}\right)\end{align*}
Consequently, we have $\lim_{n\rightarrow\infty}\sigma_{n}(\mathbb{T})=t/2$,
and $\lim_{n\rightarrow\infty}\gamma_{n}=\gamma$. In particular,
we deduce that the family $\{\sigma_{n}\}_{n=1}^{\infty}$ is tight.
Let $\sigma$ be a weak cluster point of $\{\sigma_{n}\}_{n=1}^{\infty}$,
and suppose that a subsequence $\sigma_{n_{j}}$ converges weakly
to $\sigma$ as $j\rightarrow\infty$. Then we have $\sigma(\mathbb{T})=t/2$.
Moreover, Theorems 3.5 and 4.1 yield that \[
\gamma^{p}\exp\left(-\frac{t}{2}p^{2}\right)=\widehat{\nu}(p)=\gamma^{p}\exp\left(\int_{\mathbb{T}}\frac{\zeta^{^{p}}-1-ip\Im\zeta}{1-\Re\zeta}\, d\sigma(\zeta)\right),\qquad p\in\mathbb{Z}.\]
Taking the absolute value on both sides, we have\begin{align*}
0 & =\frac{t}{2}p^{2}-\int_{\mathbb{T}}\frac{1-\Re\zeta^{^{p}}}{1-\Re\zeta}\, d\sigma(\zeta)\\
 & =\sigma(\mathbb{T})p^{2}-\int_{\mathbb{T}}\frac{1-\Re\zeta^{^{p}}}{1-\Re\zeta}\, d\sigma(\zeta)=\int_{\mathbb{T}}\left[p^{2}-\frac{1-\Re\zeta^{^{p}}}{1-\Re\zeta}\right]\, d\sigma(\zeta),\end{align*}
for every $p\in\mathbb{Z}$. Therefore, we deduce that $p^{2}=(1-\Re\zeta^{^{p}})/(1-\Re\zeta)$
for $\sigma$-almost all $\zeta\in\mathbb{T}$. Since the function
$\zeta\mapsto(1-\Re\zeta^{^{p}})/(1-\Re\zeta)$ achieves its maximum
$p^{2}$ only at $\zeta=1$, we conclude that \[
\sigma=\frac{t}{2}\delta_{1}.\]
Hence, the full sequence $\sigma_{n}$ must converge weakly to $\sigma$
because $\sigma$ is unique. The result now follows by Theorem 3.5. 
\end{proof}
\begin{rem*}
The attentive reader might have noticed that a crucial step in the
proof of Corollary 4.2 is that (4.1) uniquely determines the measure
$\sigma$. The following example inspired by \cite[Chapter IV, Section 8]{Par}
shows that this phenomenon does not happen in general. Consider the
function\[
f(\zeta)=4\pi\Im\zeta,\qquad\zeta\in\mathbb{T}.\]
Note that we have $\int_{\mathbb{T}}\zeta^{^{p}}f(\zeta)\, dm(\zeta)=2p\pi i$
when $p=\pm1$, and $\int_{\mathbb{T}}\zeta^{^{p}}f(\zeta)\, dm(\zeta)=0$
for other $p$'s. Denote by $f^{+}$ the positive part of $f$, and
by $f^{-}$ the negative part of $f$. Let us introduce measures \[
d\sigma_{1}(\zeta)=(1-\Re\zeta)f^{+}(\zeta)\, dm(\zeta),\qquad d\sigma_{2}(\zeta)=(1-\Re\zeta)f^{-}(\zeta)\, dm(\zeta).\]
Then $\sigma_{1}\neq\sigma_{2}$, and yet \[
\exp\left(\int_{\mathbb{T}}\frac{\zeta^{^{p}}-1-ip\Im\zeta}{1-\Re\zeta}\, d\sigma_{1}(\zeta)\right)=\exp\left(\int_{\mathbb{T}}\frac{\zeta^{^{p}}-1-ip\Im\zeta}{1-\Re\zeta}\, d\sigma_{2}(\zeta)\right)\]
for every $p\in\mathbb{Z}$.  
\end{rem*}
We conclude this section by showing a result concerning the weak convergence
to Haar measure $m$.

\begin{thm}
The sequence\[
\lambda_{n}\prod_{k=1}^{k_{n}}\int_{\mathbb{T}}\zeta\, d\mu_{nk}(\zeta)\]
converges to zero as $n\rightarrow\infty$ if and only if the sequence
$\mu_{n}$ converges weakly to $m$ as $n\rightarrow\infty$. 
\end{thm}
\begin{proof}
Define the measures $\sigma_{n}$ and the complex numbers $\gamma_{n}$
as in the proof of Theorem 4.1. Then Lemma 2.1 and the proof of Corollary
4.2 show that the sequence $\gamma_{n}e^{-\sigma_{n}(\mathbb{T})}$
converges if and only if the sequence \[
\widehat{\nu_{n}}(1)=\lambda_{n}\prod_{k=1}^{k_{n}}\int_{\mathbb{T}}\zeta\, d\mu_{nk}(\zeta)\]
converges. Moreover, the two sequences have the same limit. Therefore,
the result follows at once by Theorem 3.6.
\end{proof}
\begin{rem*}
Theorem 4.3 shows that if the measures $\nu_{n}$ converge weakly
to Haar measure $m$, then the measures $\mu_{n}$ converge weakly
to $m$ as well. The example below indicates that the converse of
this fact may not be true in general. Define \[
\rho_{n}=\left(1-\frac{1}{n}\right)\delta_{1}+\frac{1}{n}\delta_{-1},\qquad n\in\mathbb{N}.\]
Note that we have \[
\widehat{\rho_{n}}(p)=\begin{cases}
1 & \text{if}\; p\;\text{is even};\\
1-\frac{2}{n} & \text{if}\; p\;\text{is odd}.\end{cases}\]
Theorem 3.6 shows that the sequence $\underbrace{\rho_{n}\utimes\rho_{n}\utimes\cdots\utimes\rho_{n}}_{n^{2}\:\text{times}}$
converges weakly to $m$ as $n\rightarrow\infty$. However, the sequence
$\underbrace{\rho_{n}\circledast\rho_{n}\circledast\cdots\circledast\rho_{n}}_{n^{2}\:\text{times}}$
converges weakly to the probability measure \[
\nu=\frac{1}{2}(\delta_{1}+\delta_{-1})\]
as $n\rightarrow\infty$. Note that the limit law $\nu$ is $\circledast$-infinitely
divisible because $\nu\circledast\nu=\nu$. However, the measure $\nu$
is neither $\Utime$-infinitely divisible nor $\boxtimes$-infinitely
divisible.
\end{rem*}

\section{Measures on $\mathbb{R}$}

Let $\{\nu_{nk}:n\in\mathbb{N},1\leq k\leq k_{n}\}\subset\mathcal{M}_{\mathbb{R}}$
be an infinitesimal array. Define the probability measures $\nu_{nk}^{\circ}$
by \[
d\nu_{nk}^{\circ}(t)=d\nu_{nk}(t+a_{nk}),\]
where the numbers $a_{nk}\in[-1,1]$ are given by \[
a_{nk}=\int_{\left|t\right|<1}t\, d\nu_{nk}(t).\]
Note that the array $\{\nu_{nk}^{\circ}\}_{n,k}$ is infinitesimal,
and that $\lim_{n\rightarrow\infty}\max_{1\leq k\leq k_{n}}\left|a_{nk}\right|=0$.
We introduce analytic functions\[
f_{nk}(z)=\int_{-\infty}^{\infty}\frac{t}{1+t^{2}}\, d\nu_{nk}^{\circ}(t)+\int_{-\infty}^{\infty}\left[\frac{1+tz}{z-t}\right]\,\frac{t^{2}}{1+t^{2}}\, d\nu_{nk}^{\circ}(t),\qquad z\in\mathbb{C}^{+},\]
and note that\[
f_{nk}(z)=\int_{-\infty}^{\infty}\frac{tz}{z-t}\, d\nu_{nk}^{\circ}(t)\]
for every $n$ and $k$. Moreover, observe that $\Im f_{nk}(z)<0$
for all $z\in\mathbb{C}^{+}$ unless the measure $\nu_{nk}^{\circ}=\delta_{0}$,
and that $f_{nk}(z)=o(\left|z\right|)$ as $z\rightarrow\infty$ nontangentially.
The following result is analogous to Lemma 3.1. 

\begin{lem}
Let $\Gamma_{\alpha,\beta}$ be a truncated cone. Then for sufficiently
large $n$, we have\[
E_{\nu_{nk}^{\circ}}(z)=f_{nk}(z+a_{nk})(1+v_{nk}(z)),\]
where the sequence \[
v_{n}(z)=\max_{1\leq k\leq k_{n}}\left|v_{nk}(z)\right|\]
has properties that $\lim_{n\rightarrow\infty}v_{n}(z)=0$ for all
$z\in\Gamma_{\alpha,\beta}$, and that $v_{n}(z)=o(1)$ uniformly
in $n$ as $\left|z\right|\rightarrow\infty$, $z\in\Gamma_{\alpha,\beta}$.
\end{lem}
\begin{proof}
It was shown in \cite[Proposition 6.1]{BerPataStable} that the function
$E_{\nu_{nk}^{\circ}}(z)$ can be approximated by the function $f_{nk}(z)$
in the way we stated in the current lemma for sufficiently large $n$.
To prove the lemma, we only need to show that the function $f_{nk}(z+a_{nk})$
can be approximated by the function $f_{nk}(z)$ in the same way.
As in Lemma 3.1, we may assume that $\Im f_{nk}(z)<0$ for all $n$,
$k$, and $z\in\Gamma_{\alpha,\beta}$. Then it suffices to show that
the sequence \[
u_{n}(z)=\max_{1\leq k\leq k_{n}}\left|\frac{f_{nk}(z+a_{nk})}{f_{nk}(z)}-1\right|\]
converges to zero as $n\rightarrow\infty$ for every $z\in\Gamma_{\alpha,\beta}$,
and that $u_{n}(z)=o(1)$ uniformly in $n$ as $z\rightarrow\infty$,
$z\in\Gamma_{\alpha,\beta}$. Indeed, we have, for all $n$, $k$,
and $z\in\Gamma_{\alpha,\beta}$, that\begin{align*}
\left|f_{nk}(z+a_{nk})-f_{nk}(z)\right| & \leq\left|a_{nk}\right|\int_{-\infty}^{\infty}\frac{t^{2}}{\left|z+a_{nk}-t\right|\left|z-t\right|}\, d\nu_{nk}^{\circ}(t)\\
 & =\left|a_{nk}\right|\int_{-\infty}^{\infty}\frac{t^{2}}{\left|z-t\right|^{2}}\,\frac{\left|z-t\right|}{\left|z+a_{nk}-t\right|}\, d\nu_{nk}^{\circ}(t)\\
 & \leq2\sqrt{1+\alpha^{2}}\left|a_{nk}\right|\int_{-\infty}^{\infty}\frac{t^{2}}{\left|z-t\right|^{2}}\, d\nu_{nk}^{\circ}(t),\end{align*}
while\[
\left|f_{nk}(z)\right|\geq\left|\Im f_{nk}(z)\right|>\Im z\int_{-\infty}^{\infty}\frac{t^{2}}{\left|z-t\right|^{2}}\, d\nu_{nk}^{\circ}(t).\]
 Hence, we conclude that\[
\left|\frac{f_{nk}(z+a_{nk})}{f_{nk}(z)}-1\right|\leq\frac{\left|f_{nk}(z+a_{nk})-f_{nk}(z)\right|}{\left|\Im f_{nk}(z)\right|}\leq2\sqrt{1+\alpha^{2}}\frac{\left|a_{nk}\right|}{\Im z}.\]
 The result follows since $\lim_{n\rightarrow\infty}\max_{1\leq k\leq k_{n}}\left|a_{nk}\right|=0$.
\end{proof}
As shown in \cite[Lemma 3.1]{BerJcAdd}, the functions $f_{nk}(z)$
possess remarkable features as follows. For $y\geq1$, and for sufficiently
large $n$, we have\[
\left|\Re f_{nk}(iy)\right|\leq(3+6y)\left|\Im f_{nk}(iy)\right|,\qquad1\leq k\leq k_{n},\]
and \[
\left|\Re\left[f_{nk}(iy)-b_{nk}(y)\right]\right|\leq2\left|\Im f_{nk}(iy)\right|,\qquad1\leq k\leq k_{n},\]
where the real-valued function $b_{nk}(y)$ is given by \[
b_{nk}(y)=\int_{\left|t\right|\geq1}\left[a_{nk}+\frac{(t-a_{nk})y^{2}}{y^{2}+(t-a_{nk})^{2}}\right]\, d\nu_{nk}(t).\]

\begin{prop}
Let $\{ c_{n}\}_{n=1}^{\infty}$ be a sequence of real numbers. 
\begin{enumerate}
\item For any $y\geq1$, the sequence $\{ c_{n}+\sum_{k=1}^{k_{n}}E_{\nu_{nk}}(iy)\}_{n=1}^{\infty}$
converges if and only if the sequence $\{ c_{n}+\sum_{k=1}^{k_{n}}\left[a_{nk}+f_{nk}(iy)\right]\}_{n=1}^{\infty}$
converges. Moreover, the two sequences have the same limit.
\item If \[
L=\sup_{n\geq1}\sum_{k=1}^{k_{n}}\int_{-\infty}^{\infty}\frac{t^{2}}{1+t^{2}}\, d\nu_{nk}^{\circ}(t)<+\infty,\]
then $c_{n}+\sum_{k=1}^{k_{n}}E_{\nu_{nk}}(iy)=o(y)$ uniformly in
$n$ as $y\rightarrow\infty$ if and only if $c_{n}+\sum_{k=1}^{k_{n}}\left[a_{nk}+f_{nk}(iy)\right]=o(y)$
uniformly in $n$ as $y\rightarrow\infty$.
\end{enumerate}
\end{prop}
\begin{proof}
Fix $y\geq1$. Since $E_{\nu_{nk}^{\circ}}(z)=E_{\nu_{nk}}(z+a_{nk})-a_{nk}$,
we obtain, from Lemma 5.1, that\[
-E_{\nu_{nk}}(iy)+a_{nk}=-f_{nk}(iy)(1+u_{nk}(iy)),\]
where the sequence $u_{n}(iy)=\max_{1\leq k\leq k_{n}}\left|u_{nk}(iy)\right|$
converges to zero as $n\rightarrow\infty$. Thus, (1) follows from
(2.4) and (2.5) by setting $z_{nk}=-iE_{\nu_{nk}}(iy)+ia_{nk}$, $w_{nk}=-if_{nk}(iy)$
and $s_{nk}=0$.

Now, let us prove (2). Since $\lim_{n\rightarrow\infty}\max_{1\leq k\leq k_{n}}\left|a_{nk}\right|=0$
and $u_{n}(iy)=o(1)$ uniformly in $n$ as $y\rightarrow\infty$,
we may assume that $\left|a_{nk}\right|\leq1/2$, and that $u_{n}(iy)<1/6$,
for all $n$, $k$ and for sufficiently large $y$. Observe that \begin{eqnarray*}
\sum_{k=1}^{k_{n}}\left|b_{nk}(y)\right| & = & \sum_{k=1}^{k_{n}}\left|\int_{\left|t\right|\geq1}\left[a_{nk}+\frac{(t-a_{nk})y^{2}}{y^{2}+(t-a_{nk})^{2}}\right]\, d\nu_{nk}(t)\right|\\
 & \leq & (1+y)\sum_{k=1}^{k_{n}}\int_{\left|t\right|\geq1}\frac{1}{2}\, d\nu_{nk}(t)\leq5y\sum_{k=1}^{k_{n}}\int_{\left|t\right|\geq1}\frac{1}{5}\, d\nu_{nk}(t)\\
 & \leq & 5y\sum_{k=1}^{k_{n}}\int_{\left|t\right|\geq1}\frac{(t-a_{nk})^{2}}{1+(t-a_{nk})^{2}}\, d\nu_{nk}(t)\leq5yL.\end{eqnarray*}
Then (2.4) and (2.5) imply that \[
\left|\left(\sum_{k=1}^{k_{n}}E_{\nu_{nk}}(iy)\right)-\left(\sum_{k=1}^{k_{n}}\left[a_{nk}+f_{nk}(iy)\right]\right)\right|\leq\frac{1}{2}\left|\sum_{k=1}^{k_{n}}\Im f_{nk}(iy)\right|+5yLu_{n}(iy),\]
and \[
\frac{1}{2}\left|\sum_{k=1}^{k_{n}}\Im f_{nk}(iy)\right|\leq\left|\sum_{k=1}^{k_{n}}\Im E_{\nu_{nk}}(iy)\right|+5yLu_{n}(iy),\]
for $n\in\mathbb{N}$. Then (2) follows since $u_{n}(iy)=o(1)$ uniformly
in $n$ as $y\rightarrow\infty$.
\end{proof}
We are now ready for the main result of this section. With Proposition
5.2 in hands, one can follow almost word for word the argument of
\cite[Theorem 3.3]{BerJcAdd} to prove the following result. Therefore,
we will not repeat this rather lengthy proof here but refer to \cite{BerJcAdd}
for its details. 

\begin{thm}
Fix a real number $\gamma$ and a finite positive Borel measure $\sigma$
on $\mathbb{R}$. Let $\{ c_{n}\}_{n=1}^{\infty}$ be a sequence of
real numbers. Then the following statements are equivalent:
\begin{enumerate}
\item The sequence $\delta_{c_{n}}*\nu_{n1}*\nu_{n2}*\cdots*\nu_{nk_{n}}$
converges weakly to $\nu_{*}^{\gamma,\sigma}$.
\item The sequence $\delta_{c_{n}}\boxplus\nu_{n1}\boxplus\nu_{n2}\boxplus\cdots\boxplus\nu_{nk_{n}}$
converges weakly to $\nu_{\boxplus}^{\gamma,\sigma}$.
\item The sequence $\delta_{c_{n}}\uplus\nu_{n1}\uplus\nu_{n2}\uplus\cdots\uplus\nu_{nk_{n}}$
converges weakly to $\nu_{\uplus}^{\gamma,\sigma}$.
\item The sequence of measures\[
d\sigma_{n}(t)=\sum_{k=1}^{k_{n}}\frac{t^{2}}{1+t^{2}}\, d\nu_{nk}^{\circ}(t)\]
 converges weakly on $\mathbb{R}$ to $\sigma$, and the sequence
of numbers\[
\gamma_{n}=c_{n}+\sum_{k=1}^{k_{n}}\left[a_{nk}+\int_{-\infty}^{\infty}\frac{t}{1+t^{2}}\, d\nu_{nk}^{\circ}(t)\right]\]
converges to $\gamma$ as $n\rightarrow\infty$. 
\end{enumerate}
\end{thm}

\section*{acknowledgments}

The author would like to thank his adviser Professor Hari Bercovici
for encouragement and many helpful discussions during his doctoral
studies at Indiana University.

\address{\textsc{\scriptsize Department of Mathematics, Indiana University
, Bloomington, Indiana 47405 , U.S.A.}}

\email{jiuwang@indiana.edu}
\end{document}